\documentclass[12pt,a4paper]{article}
\usepackage[english]{babel}
\usepackage{amssymb}
\usepackage{amsmath}
\usepackage{amsfonts}
\usepackage{bbm}
\usepackage{textcomp}
\usepackage{color}

\def\theequation{\arabic{section}.\arabic{equation}}

\newcommand{\be}{\begin{equation}}
\newcommand{\ee}{\end{equation}}
\newcommand{\bd}{\begin{displaymath}}
\newcommand{\ed}{\end{displaymath}}
\newcommand{\ba}{\begin{eqnarray}}
\newcommand{\ea}{\end{eqnarray}}

\newcommand{\p}{\partial}

\def\N{\mathbb{N}}

\def\R{\mathbb{R}}

\def\sign{\,{\rm sign}\, }

\def\a{\alpha}
\def\b{\beta}
\def\d{\delta}
\def\D{\Delta}
\def\e{\varepsilon}
\def\f{\varphi}
\def\g{\gamma}

\def\l{\lambda}

\def\U{\Upsilon}

\def\calC{\mathcal{C}}
\def\calF{\mathcal{F}}

\def\calP{\mathcal{P}}
\def\calT{\mathcal{T}}

\def\rmd{\mathrm{d}}
\def\rme{\mathrm{e}}
\def\rmL{\mathrm{L}}

\def\div{\mathrm{div}}

\newcommand{\email}[1]{{\small E-mail: {\textsf {#1}}}}
\newcommand{\http}[1]{{\small Internet: {\textsf {#1}}}}

\newtheorem{theorem}{Theorem}[section]
\newtheorem{lemma}[theorem]{Lemma}
\newtheorem{proposition}[theorem]{Proposition}
\newtheorem{corollary}[theorem]{Corollary}
\newtheorem{remark}[theorem]{Remark}

\def\proof{\par{\it \noindent Proof.-} \ignorespaces}
\def\endproof{\hfill{{\ \vbox{\hrule\hbox{%
   \vrule height1.3ex\hskip0.8ex\vrule}\hrule
  }}\par}}

\begin{document}

\author{Jos\'e A. Carrillo\footnote{ICREA (Instituci\'o Catalana de
Recerca i Estudis Avan\c cats) and Departament de
Mate\-m\`a\-tiques, Universitat Aut\`onoma de Barcelona, E--08193
Bellaterra, Spain. \email{carrillo@mat.uab.cat},
\http{http://kinetic.mat.uab.es/$\sim$carrillo/}} \kern8pt
Philippe Lauren\c{c}ot \footnote{Institut de
Math\'ematiques de Toulouse, CNRS (UMR~5219) \& Universit\'e de
Toulouse, 118 route de Narbonne, F--31062 Toulouse C\'edex 9,
France. \email{laurenco@mip.ups-tlse.fr},
\http{http://www.mip.ups-tlse.fr/$\sim$laurenco/}} \kern8pt\&
\kern8pt Jes\'us Rosado\footnote{Departament de
Mate\-m\`a\-tiques, Universitat Aut\`onoma de Barcelona, E--08193
Bellaterra, Spain. \email{jrosado@mat.uab.cat}}}
\date{\today}

\title{Fermi-Dirac-Fokker-Planck Equation: Well-posedness \& Long-time
Asymptotics} 

\maketitle

\begin{abstract}
A Fokker-Planck type equation for interacting particles with
exclusion principle is analysed. The nonlinear drift gives rise to
mathematical difficulties in controlling moments of the
distribution function. Assuming enough initial moments are finite,
we can show the global existence of weak solutions for this
problem. The natural associated entropy of the equation is the
main tool to derive uniform in time a priori estimates for the
kinetic energy and entropy. As a consequence, long-time
asymptotics in $L^1$ are characterized by the Fermi-Dirac
equilibrium with the same initial mass. This result is achieved
without rate for any constructed global solution and with
exponential rate due to entropy/entropy-dissipation arguments for
initial data controlled by Fermi-Dirac distributions. Finally,
initial data below radial solutions with suitable decay at
infinity lead to solutions for which the relative entropy towards
the Fermi-Dirac equilibrium is shown to converge to zero without
decay rate.
\end{abstract}

\section{Introduction}
Kinetic equations for interacting particles with exclusion
principle, such as fermions, have been introduced in the physics
literature in \cite{bib:f1,bib:k,bib:k2,bib:k3,bib:lkq,bib:rk} and
the review \cite{bib:f2}. Spatially inhomogeneous equations appear
from formal derivations of generalized Boltzmann equations and
Uehling-Uhlenbeck kinetic equations both for fermionic and bosonic
particles. The most relevant questions related to these problems
concern their long-time asymptotics and the rate of convergence
towards global equilibrium if any.

The spatially inhomogeneous situation has been recently studied in
\cite{bib:ns2}, where the long time asympotics
of these models in the torus is shown to be given by spatially
homogeneous equilibrium given by Fermi-Dirac distributions when the
initial data is not far from equilibrium in a suitable Sobolev
space. This nice result is based on techniques developed in previous
works \cite{bib:mn,bib:ns}. Other related mathematical results for
Boltzmann-type models have appeared in \cite{bib:d,bib:lw}. 

In this work, we focus on the global existence of solutions and
the convergence of solutions towards global equilibrium in the
spatially homogeneous case without any smallness assumption on the
initial data. Preliminary results in the one-dimensional setting
were reported in \cite{bib:crs}. More precisely, we analyse in
detail the following Fokker-Planck equation for fermions, see for
instance \cite{bib:f2},
\begin{equation}\label{eq:fpf}
    \frac{\p f}{\p t}=\Delta_v f + {\rm div}_v[vf(1-f)], \qquad v \in \R^N, t>0,
\end{equation}
with initial condition $f(0,v)=f_0(v)\in L^1(\R^N)$, $0\leq
f_0\leq 1$ and suitable moment conditions to be specified below.
Here, $f=f(t,v)$ is the density of particles with velocity $v$ at
time $t\geq 0$.

This equation has been proposed in order to describe the dynamics
of classical interacting particles, obeying the
exclusion-inclusion principle in \cite{bib:k}. In fact, equation
\eqref{eq:fpf} is formally equivalent to
$$
\frac{\p f}{\p
t}=\div_v\left[f(1-f)\nabla_v\left(\log\left(\frac{f}{1-f}\right)+\frac{|v|^2}{2}\right)\right]
$$
from which it is easily seen that Fermi-Dirac distributions defined by

$$
F^{\beta}(v):=\frac{1}{1+\beta\mathrm{e}^{\frac{|v|^2}{2}}}
$$
with $\beta\geq 0$ are stationary solutions. Moreover, for each value of $M\ge 0$, there exists a unique $\beta=\beta(M) \geq 0$
such that $F^{\beta(M)}$ has mass $M$, that is, $\|F^{\beta(M)}\|_1=M$. Throughout the paper we shall denote $F^{\beta(M)}$ by $F_M$.

Another striking property of this equation is the existence of a
formal Liapunov functional, related to the standard entropy
functional for linear and nonlinear Fokker-Planck models
\cite{bib:ct,bib:cjmtu}, given by
$$
H(f):=\frac{1}{2}\int_{\R^N}|v|^2f(v)\,\rmd v+\int_{\R^N}
[(1-f)\log(1-f)+f\log(f)]\,\rmd v.
$$
We will show that this functional plays the same role as the
H-functional for the spatially homogeneous Boltzmann equation, see for instance \cite{bib:t}. In particular it will be crucial to characterize long-time asymptotics of \eqref{eq:fpf}. In fact, the entropy method
will be the basis of the main results in this work; more precisely by taking the formal time derivative of $H(f)$, we conclude that
$$
\frac{d}{dt}H(f)=-\int_{\R^N}
f(1-f)\left|v+\nabla_v
\log\left(\frac{f}{1-f}\right)\right|^2\,\rmd v\leq 0.
$$
Therefore, to show the global equilibration of solutions to \eqref{eq:fpf} we need to find the right functional setting to show the
entropy dissipation. Furthermore, if we succeed in relating
functionally the entropy and the entropy dissipation, we will be
able to give decay rates towards equilibrium. These are the
main objectives of this work. Let us finally mention that these
equations are of interest as typical examples of gradient flows
with respect to euclidean Wasserstein distance of entropy
functionals with nonlinear mobility, see \cite{bib:bfd,bib:cmv}
for other examples and related problems.

In section 2, we will show the global existence of solutions for
equation \eqref{eq:fpf} based on fixed point arguments, estimates
involving moment bounds and the conservation of certain properties
of the solutions. The suitable functional setting is reminiscent of
the one used in equations sharing a similar structure and
technical difficulties as those treated in \cite{bib:ez,bib:gw}.
The main technical obstacle for the Fermi-Dirac-Fokker-Planck
equation \eqref{eq:fpf} lies in the control of moments. Next, in
section 3, we show that the constructed solutions verify that
the entropy is decreasing, and from that, we prove the convergence
towards global equilibrium without rate. Again, here the
uniform-in-time control of the second moment is crucial. Finally,
we obtain an exponential rate of convergence towards equilibrium
if the initial data is controlled by Fermi-Dirac distributions
and the convergence to zero of the relative entropy when controlled by radial solutions.

\section{Global Existence of Solutions}\label{s:EoS}
\setcounter{equation}{0}

In this section, we will show the global existence of solutions
to the Cauchy problem to \eqref{eq:fpf}. We start by proving
local existence of solutions together with a characterization of
the time-span of these solutions. Later, we show further
regularity properties of these solutions with the help of estimates on
derivatives. Based on these estimates we can derive further
properties of the solutions: conservation of mass, positivity,
$L^\infty$ bounds, comparison principle, moment estimates and
entropy estimates. All of these uniform estimates allow us to show
that solutions can be extended and thus exist for all times.

\subsection{Local Existence}\label{s:lEiU}

We will prove the local existence and uniqueness of solutions using
contraction-principle arguments as in \cite{bib:bfd,bib:ez,bib:gw} for
instance. As a first step, let us note that we can write
\eqref{eq:fpf} as
\begin{equation}\label{eq:fpfFon}
\frac{\p f}{\p t}=\D_v f+\div_v(vf)-\div_v(vf^2)
\end{equation}
and, due to Duhamel's formula, we are led to consider the
corresponding integral equation
\begin{equation}\label{eq:fpfFonInt}
f(t,v)=\int_{\R^N}\calF(t,v,w)f_0(w)\rmd w -\int_0^t\int_{\R^N}\calF(t-s,v,w)(\div_w(wf(s,w)^2))\,\rmd w\,\rmd s
\end{equation}
where $\calF(t,v,w)$ is the fundamental solution for the
homogeneous Fokker-Planck equation:
$$
\frac{\p f}{\p t}=\div_v(vf+\nabla_vf)
$$
given by
$$
\calF(t,v,w):=a(t)^{-\frac{N}{2}}M_{\nu(t)}(a(t)^{-\frac{1}{2}}v-w)
$$
with
$$
a(t):=\rme^{-2t}\quad\text{,}\quad\nu(t):=\rme^{2t}-1\quad\text{and}\quad
M_{\l}(\xi):=(2\pi \l)^{-\frac{N}{2}}\rme^{-\frac{|\xi|^2}{2\l}}
$$
for any $\lambda>0$. Let us define the operator $\calF(t,v)[g]$
acting on functions $g$ as:
\begin{equation}\label{eq:fpFonSol}
\calF(t,v)[g(w)]= \int_{\R^N}\calF(t,v,w)g(w)\,\rmd w \text{.}
\end{equation}

Note that by integration by parts, the expression
$\calF(t,v)[\div_w(wf^2(w))]$ is equivalent to:
\begin{align*}
\int_{\R^N}\Bigg{(}\frac{\rme^{N t}}
{\left(2\pi\left(\rme^{2t}-1\right)\right)^{\frac{N}{2}}}
\,\rme^{-\frac{|\rme^{t}v-w|^2} {2\left(\rme^{2t}-1\right)}}
\Bigg{)}&\div_w(wf(w)^2)\,\rmd w\\
&=-\int_{\R^N}\left[\nabla_w\left(
\frac{\rme^{N t}} {\left(2\pi\left(\rme^{2t}-1\right)\right)^{\frac{N}{2}}}
\,\rme^{-\frac{|\rme^{t}v-w|^2}{2\left(\rme^{2t}-1\right)}}
\right)\cdot w\right] f(w)^2\,\rmd w\\
&=-\int_{\R^N}\rme^{-t}\left( \nabla_v \calF(t,v,w)\cdot
w\right)f(w)^2\,\rmd w\\
&=:-\rme^{-t}\,\nabla_v\calF(t,v)[wf(w)^2]
\end{align*}
so that (\ref{eq:fpfFonInt}) becomes
\begin{equation}\label{eq:fpfFonInt2}
f(t,v)=\calF(t,v)[f_0(w)]
+\int_0^t\rme^{-(t-s)}\nabla_v\calF(t-s,v)[wf(s,w)^2]\,\rmd s\text{.}
\end{equation}

We will now define a space in which the functional induced by
(\ref{eq:fpfFonInt2})
\begin{equation}\label{eq:fpfFonIntrec}
\calT[f](t,v):=\calF(t,v)[f_0(w)] +
\int_0^t\rme^{-(t-s)}\nabla_v\calF(t-s,v)[wf(s,w)^2]\,\rmd
s
\end{equation}
has a fixed point. To this end, we define the spaces
$\Upsilon:=\rmL^{\infty}(\R^N)\cap\rmL^1_1(\R^N)\cap\rmL^p_m(\R^N)$
and $\Upsilon_T:=\calC([0,T];\Upsilon)$
with norms
$$
\|f(t)\|_{\Upsilon}:=\max\{\|f(t)\|_{\infty},\|f(t)\|_{\rmL^1_1},
\|f(t)\|_{\rmL_m^p}\} \quad\text{and} \quad
\|f\|_{\Upsilon_T}:=\max_{0\leq t\leq T}\|f(t)\|_{\Upsilon}
$$
for any $T>0$, where we omit the N-dimensional euclidean space
$\R^N$ for notational convenience and
$$
\|f\|_{\rmL_m^p}:=\|(1+|v|^m)f\|_p\quad \text{and}
\quad\|f\|_p:=\left(\int_{\R^N}|f|^p\rmd
v\right)^{\frac{1}{p}}\text{.}
$$
In the following, we will see that for $p>N$, $p\ge 2$, and $m\geq 1$ we can choose $q$ and $r$ satisfying
\begin{equation}
\frac{Np}{N+p}<\frac{p}{2}\le r \le \frac{mp}{m+1} < p \qquad \mbox{and} \qquad \frac{p}{2} \leq q\leq p \label{eq:conditions}
\end{equation}
such that
$\|\calT[f]\|_{\Upsilon_T}$ is bounded by $\|f\|_{\Upsilon_T}$.
Let us fix such parameters $p$, $m$, $r$, $q$ and $0\leq t\leq T$. Due to Proposition
\ref{prop:fpnucpqbounds} and $q\leq p\leq
2q$, we can compute
\begin{align*}
\|\calT[f](t)\|_{\infty}&\leq C\rme^{Nt}\|f_0\|_{\infty} + \int_0^tC\frac{\rme^{N(t-s)}}{\nu(t-s)^{\frac{N}{2q}+\frac{1}{2}}}\,\||w|f^2(s)\|_q\,\rmd s\\
&\leq C\rme^{Nt}\|f_0\|_{\infty} + \int_0^tC\frac{\rme^{N(t-s)}}{\nu(t-s)^{\frac{N}{2q}+\frac{1}{2}}}\, \|f(s)\|_{\infty}^{2-\frac{p}{q}}\|f(s)\|_{\rmL^p_m} ^{\frac{p}{q}}\rmd s\\
&\leq C\rme^{Nt}\|f_0\|_{\infty} +
\int_0^tC\frac{\rme^{N(t-s)}}{\nu(t-s)^{\frac{N}{2q}+\frac{1}{2}}}
\,\rmd s\, \|f\|_{\U_T}^2\\
&\leq C\rme^{Nt}\|f_0\|_{\infty} + C\, \mathcal{I}_ 1(t)\, \|f\|_{\U_T}^2\text{,}
\end{align*}
where
$$
\mathcal{I}_1(t):=\int_{\rme^{-2t}}^1
\chi^{-\frac{1}{2}\left(N-\frac{N}{q}-1\right)-1}(1-\chi)^{-\frac{1}{2}(\frac{N}{q}+1)}
\rmd\chi<\infty
$$
by the choice \eqref{eq:conditions} of $q$. In the same way, since $r$ satisfies $(m+1)r\leq mp$ and
$2r\geq p$, we get
\begin{align*}
\|\calT[f](t)\|_{\rmL_m^p}&\leq C\rme^{\frac{N}{p'}t}\|f_0\|_{\rmL^p_m} + \int_0^t C \frac{\rme^{\frac{N}{p'}(t-s)}} {\nu(t-s)^{\frac{N}{2}\left(\frac{1}{r} - \frac{1}{p}\right) + \frac{1}{2}}}\||w|f^2(s)\|_{\rmL^r_m}\rmd s\\
&\leq C\rme^{\frac{N}{p'}t}\|f_0\|_{\rmL^p_m} + \int_0^t C \frac{\rme^{\frac{N}{p'}(t-s)}} {\nu(t-s)^{\frac{N}{2}\left(\frac{1}{r} -\frac{1}{p}\right)+\frac{1}{2}}}\|f(s)\|_{\infty}^{2-\frac{p}{r}} \|f(s)\|_{\rmL^p_m}^{\frac{p}{r}} \rmd s\\
&\leq C\rme^{\frac{N}{p'}t}\|f_0\|_{\rmL^p_m} +
\int_0^tC\frac{\rme^{\frac{N}{p'}(t-s)}}{\nu(t-s)^{\frac{N}{2}\left(\frac{1}{r}-\frac{1}{p}\right)+\frac{1}{2}}}\,
\rmd s \,\|f\|_{\U_T}^2 \\
&\leq C\rme^{\frac{N}{p'}t}\|f_0\|_{\rmL^p_m} + C\, \mathcal{I}_2(t)\, \|f\|_{\U_T}^2\text{,}
\end{align*}
where
$$
\mathcal{I}_2(t) := \int_{\rme^{-2t}}^1
\chi^{-\frac{1}{2}\left[\frac{N}{p'}-\left(N\left(\frac{1}{r}-
\frac{1}{p}\right)+1\right)\right]-1}(1-\chi)^{-\frac{N}{2}\left(\frac{1}{r}-\frac{1}{p}\right)-\frac{1}{2}} \rmd\chi<\infty
$$
by the choice \eqref{eq:conditions} of $r$.

Finally we can estimate
\begin{align*}
\|\calT[f](t)\|_{\rmL^1_1}\leq
C\|f_0\|_{\rmL^1_1} +
\int_0^t\frac{C}{\nu(t-s)^{\frac{1}{2}}}\||w|f^2(s)\|_{\rmL^1_1}
\,\rmd s
\end{align*}
where by interpolation, we get as $p\geq 2$ and $m\geq 1$
\begin{align}\label{eq:1norm_wf2}
\||w|f^2\|_{\rmL^1_1}&=\int_{\R^N} (1+|w|)|w|f^2 \,\rmd w
\leq\int_{\R^N}(1+|w|)^2f^2\,\rmd w\nonumber\\&\leq
\left(\int_{\R^N}(1+|w|)f\,\rmd w\right)^{\frac{p-2}{p-1}}
\left(\int_{\R^N}(1+|w|)^pf^p\,\rmd w\right)^{\frac{1}{p-1}}\nonumber\\
&\leq
\|f\|_{\rmL^1_1}^{\frac{p-2}{p-1}}\|f\|_{\rmL_m^p}^{\frac{p}{p-1}}.
\end{align}

{\noindent Consequently}
\begin{equation*}
\|\calT[f](t)\|_{\rmL^1_1}\leq C\|f_0\|_{L^1_1}+C\int_{\rme^{-2t}}^1
\chi^{-\frac{3}{2}}(1-\chi)^{-\frac{1}{2}}\rmd\chi\,\|f\|_{\Upsilon_T}^2.
\end{equation*}
We next check the existence of a fixed point of
(\ref{eq:fpfFonIntrec}) in $\Upsilon_T$. To this end, we define a
sequence $(f_n)_{n\ge 1}$ by $f_{n+1}=\calT[f_n]$ for $n\ge
0$. Collecting all the above estimates, we can write
$$
\|f_{n+1}(t)\|_{\Upsilon}\leq
C_1(N,t)\|f_0\|_{\Upsilon}+C_2(N,p,q,r,t)\|f_n\|_{\Upsilon_T}^2
$$
for any $0\leq t\leq T$ and any $T>0$, with
$$
\begin{array}{rcll}
C_1(N,t)&:=&C\rme^{Nt}\\
C_2(N,p,q,r,t)&:=&C\max\Bigg{\{}\!\!\!&\displaystyle{\mathcal{I}_1(t)
, \mathcal{I}_2(t)},\displaystyle{\int_{\rme^{-2t}}^1
\chi^{-\frac{3}{2}}(1-\chi)^{-\frac{1}{2}} \rmd\chi}\Bigg{\}}
\end{array}$$
which are clearly increasing with $t$ and $C_2(t)$ tends to $0$ as $t$
does. Thus, for any $T>0$
$$
\|f_{n+1}\|_{\Upsilon_T}\leq
C_1(T)\,\|f_0\|_{\Upsilon}+C_2(T)\,\|f_n\|_{\Upsilon_T}^2
$$
with $C_1(T)=C_1(N,T)$ and $C_2(T)=C_2(N,p,q,r,T)$, both being
increasing functions of $T$. We may also assume that $C_1(T)\geq 1$
without loss of generality.

{F}rom now on, we will follow the arguments in \cite{bib:lr}. We
will first show that if $T$ is small enough, the functional $\calT$ is
bounded in $\U_T$, which will in turn imply the
convergence. Let us take $T>0$ and $\d>0$ which verify
$$
\|f_0\|_{\U}<\d \qquad \mbox{and} \qquad
0<\d<\frac{1}{4C_1(T)C_2(T)}.
$$
Then, let us prove by induction that
$\|f_n\|_{\U_T}<2C_1(T)\d$ for all $n$. By the choice of $T$ and $\d$ we
have $\|f_0\|_{\U}<C_1(T)\d<2C_1(T)\d$. If we suppose that
$\|f_n\|_{\Upsilon_T}<2C_1(T)\d$, we have
$$
\|f_{n+1}\|_{\Upsilon_T}<C_1(T)\d+4C_1^2(T)C_2(T)\d^2<2C_1(T)\d,
$$
hence the claim. Now, computing the difference between two consecutive iterations of the functional and proceeding with the same estimates as above, we can see for any $0\leq t\leq T$ that
\begin{align*}
\|f_{n+1}-f_{n}\|_{\U_T}
&=\Bigg{\|}\int_0^t\rme^{-(t-s)}\nabla_v\calF(t-s,v)\left[w\left[f_n^2-f_{n-1}^2\right]\right]
\rmd s\Bigg{\|}_{\U_T}\\
&\leq C_2(T)\sup_{[0,T]}\big{\|}f_n+f_{n-1}\big{\|}_{\infty} \big{\|}f_n-f_{n-1}\big{\|}_{\U_T}\\
&\leq C_2(T)\left(\big{\|}f_n\big{\|}_{\U_T} +
\big{\|}f_{n-1}\big{\|}_{\U_T}\right)
\big{\|}f_n-f_{n-1}\big{\|}_{\U_T}\\
&\leq 4C_1(T)C_2(T)\d\big{\|}f_n-f_{n-1}\big{\|}_{\U_T}\leq
(4C_1(T)C_2(T)\d)^n\big{\|}f_1-f_0\big{\|}_{\U_T}\ .
\end{align*}
Since $4C_1(T)C_2(T)\d<1$ we can conclude that
there exists a function $f_*$ in $\U_T$ which is a fixed point for
$\calT$, and hence a solution to the integral equation
(\ref{eq:fpfFonInt}). It is not difficult to check that the solution
$f\in\U_T$ to the integral
equation is a solution of \eqref{eq:fpf} in the sense of
distributions defining our concept of solution. We summarize the
results of this subsection in the following result.

\begin{theorem}[Local Existence]\label{t:localE}
Let $m\geq 1$, $p>N$, $p\geq 2$, and $f_0\in \U$. Then there exists
$T>0$ depending only on the norm of the initial condition $f_0$ in
$\U$, such that \eqref{eq:fpf} has a unique solution $f$ in
$\calC([0,T];\U)$ with $f(0)=f_0$.
\end{theorem}

\begin{remark}
The previous theorem is also valid for $f_0\in
(\rmL^{\infty}\cap\rmL^p_m\cap\rmL^1)(\R^N)$, with a solution
defined in
$\calC([0,T];(\rmL^{\infty}\cap\rmL^p_m\cap\rmL^1)(\R^N))$ but we
will need to have the first moment of the solution bounded in
order to be able to extend it to a global in time solution. We thus
include here this additional condition.
\end{remark}

\begin{remark}
With the same arguments used to prove Theorem {\rm \ref{t:localE}}
we can prove an equi\-va\-lent result for the
Bose-Einstein-Fokker-Planck equation
\begin{equation*}
    \frac{\p f}{\p t}=\Delta_v f + {\rm div}_v[vf(1+f)], \qquad v \in \R^N, t>0.
\end{equation*}
\end{remark}

%
%
\subsection{Estimates on Derivatives}\label{s:D}
Let us now work on estimates on the derivatives. By taking the
gradient in the integral equation, we obtain
\begin{equation}
\nabla_v f(t,v)=\nabla_v\calF(t,v)[f(w)] -\int_0^t\nabla_v\calF(t-s,v)[\div_w(wf^2(s,w))]\,\rmd s.
\end{equation}
where $\nabla_v\calF(t,v)[g]$ is defined as the vector:
$$\nabla_v\calF(t,v)[g]:=\int_{\R^N}\nabla_v\calF(t,v,w)g(w)\rmd w$$
for the real-valued function $g$.
Here, we will consider a space $X_T$ with suitable weighted norms
for the derivatives
$$
\|f\|_{X_T}=\max\left\{ \|f\|_{\U_T},\sup_{0<t<
T}\nu(t)^{\frac{1}{2}}
 \|\nabla_v f(t)\|_{\rmL_m^p},\sup_{0<t<
T}\nu(t)^{\frac{1}{2}}
 \|\nabla_v f(t)\|_{\rmL^1_1}
\right\}
$$
where for notational simplicity we refer to $\||\nabla_v
f|\|_{\rmL^p_m}$ as $\|\nabla_v f\|_{\rmL^p_m}$. Let us estimate
the $\rmL^p_m$- and $\rmL^1$-norms of $\nabla_v f$ using again the
results in Proposition \ref{prop:fpnucpqbounds} as follows: for $r\in[1,p)$ satisfying \eqref{eq:conditions}
\begin{eqnarray*}
\|\nabla_v f(t)\|_{\rmL_m^p}&\leq&
C\frac{\rme^{\left(\frac{N}{p'}+1\right)t}}{\nu(t)^{\frac{1}{2}}}
\|f_0\|_{\rmL_m^p} +\int_0^t\|\nabla_v
\calF[2f(w\cdot\nabla_wf)] + Nf^2\|_{\rmL_m^p} \rmd s\\
&\leq&
C\frac{\rme^{\left(\frac{N}{p'}+1\right)t}}{\nu(t)^{\frac{1}{2}}}
\|f_0\|_{\rmL_m^p}
+C\int_0^t\frac{\rme^{\left(\frac{N}{p'}+1\right)(t-s)}}{\nu(t-s)^{\frac{1}{2}}}\|f(s)\|_{\rmL_m^p}\|f(s)\|_{\infty}
\rmd s\\& &
+C\int_0^t\frac{\rme^{\left(\frac{N}{p'}+1\right)(t-s)}}{\nu(t-s)^{\frac{N}{2}\left(\frac{1}{r}-\frac{1}{p}\right)+\frac{1}{2}}}
\|f(w\cdot\nabla_w f)\|_{\rmL_m^r}\rmd s\\
&\leq&
C\frac{\rme^{\left(\frac{N}{p'}+1\right)t}}{\nu(t)^{\frac{1}{2}}}
\|f_0\|_{\rmL_m^p}+C\|f\|_{\U_T}^2\int_{\rme^{-2t}}^1\chi^{-\frac{N+2p'}{2p'}}(1-\chi)^{-\frac{1}{2}}
\rmd s\\ & & +C\sup_{0<s<T}\left\{\nu(s)^{1/2}\|f(s)(w\cdot\nabla_w
f(s))\|_{\rmL_m^r}\right\} I(t)
\end{eqnarray*}
where
\begin{align*}
I(t)&\leq\frac{\rme^{-t}}{2}\int_{\rme^{-2t}}^1\rme^{t\left(
  \frac{N+2r'}{r'}\right)} (1-\chi)^{-(\frac{N}{2}(\frac{1}{r}-\frac{1}{p})+\frac{1}{2})} (\chi-e^{-2t})^{-\frac{1}{2}}\rmd\chi\\
&\leq\frac{1}{2}\rme^{t(\frac{N+r'}{r'})}\left[\int_{\rme^{-2t}}^{\frac{1+\rme^{-2t}}{2}}
  \left(\frac{1-\rme^{-2t}}{2}\right)^{-\frac{N}{2}(\frac{1}{r}-\frac{1}{p})-\frac{1}{2}}(\chi-\rme^{-2t})^{-\frac{1}{2}}\rmd\chi\right.\\
&\qquad\qquad\left.+\int_{\frac{1+\rme^{-2t}}{2}}^1
  (\chi-\rme^{-2t})^{-(\frac{N}{2}(\frac{1}{r}-\frac{1}{p})-\frac{1}{2}}
\left(\frac{1-\rme^{-2t}}{2} \right)^{-\frac{1}{2}}\rmd\chi \right]\\
&\leq C\rme^{t\frac{N+r'}{r'}}(1-\rme^{-2t})^{-\frac{N}{2}(\frac{1}{r}-\frac{1}{p})}\\
&\leq C\rme^{t(N-\frac{N}{r}+1+\frac{N}{r}-\frac{N}{p})}\nu(t)^{-\frac{1}{2}} \nu(t)^{\frac{1}{2}-\frac{N}{2}(\frac{1}{r}-\frac{1}{p})}\\
&\leq Ch(t)\nu(t)^{-\frac{1}{2}}
\end{align*}
with
$h(t):=\rme^{t(\frac{N+p'}{p'})}\nu(t)^{\frac{1}{2}-\frac{N}{2}(\frac{1}{r}
  -\frac{1}{p})}$ which is an increasing function of time with
$h(0)=0$ since $p>r>Np/(N+p)$.
It remains to estimate $\|f(w\cdot\nabla_w
f)\|_{\rmL_m^r}$:
$$
\|f(w\cdot\nabla_w f)\|_{\rmL_m^r} \leq C\left(\int_{\R^N}
f^r|\nabla_w f|^r\rmd w+\int_{\R^N} |w|^{(m+1)r}f^r|\nabla_w
f|^r\rmd w\right)^{\frac{1}{r}}
$$
Now, we can bound these integrals by using H\"{o}lder's inequality to obtain
$$
\int_{\R^N} f^r |\nabla_w f|^r\rmd w \leq \left(\int_{\R^N}
f^{\frac{pr}{p-r}}\rmd w\right)^{\frac{p-r}{p}}\left(\int_{\R^N}
|\nabla_w f|^p\rmd w\right)^{\frac{r}{p}}
$$
and
$$
\int_{\R^N} |w|^{(m+1)r}f^r|\nabla_w f|^r\rmd w \leq
\left(\int_{\R^N} |w|^{\frac{pr}{p-r}}f^{\frac{pr}{p-r}}\rmd w
\right)^{\frac{p-r}{p}}\left(\int_{\R^N} |w|^{mp}|\nabla_w f|^p
\rmd w \right)^{\frac{r}{p}}.
$$

Since $p<pr/(p-r)\leq mp$ or equivalently
$(m+1)r/m\leq p < 2r$ by \eqref{eq:conditions}, we have for any $0<t\leq T$
\begin{equation*}
\int_{\R^N} f^r|\nabla_w f|^r\rmd w
\leq\|f\|_{\infty}^{2r-p}\|f\|_p^{p-r}\|\nabla_w f\|_p^{r}\leq
\frac{\|f\|_{X_T}^{2r}}{\nu(t)^{\frac{r}{2}}}
\end{equation*}
and
\begin{equation*}
\int_{\R^N} |w|^{(m+1)r}f^r|\nabla_w f|^r\rmd w
\leq\|f\|_{\infty}^{2r-p}\|f\|_{\rmL_m^p}^{p-r}\|\nabla_w
f\|_{\rmL_m^p}^{r}\leq
\frac{\|f\|_{X_T}^{2r}}{\nu(t)^{\frac{r}{2}}} .
\end{equation*}
Putting together the above estimates we have shown that,
$$
\nu(t)^{1/2}\|f(t)(w\cdot\nabla_wf(t))\|_{\rmL_m^r}\leq C\|f\|_{X_T}^2
$$
and
\begin{equation}
\nu(t)^{\frac{1}{2}}\|\nabla_v f(t)\|_{\rmL_m^p}\leq
C_1^1(T,N,p)\|f_0\|_{\rmL_m^p}+C_2^1(T,N,p,r)\|f\|_{X_T}^2
\label{eq:der1}
\end{equation}
with $C_1^1$ and $C_2^1$ increasing functions of $T$ and
for any $0< t\leq T$. Analogously, we reckon
\begin{align*}
\|\nabla_v f(t)\|_{\rmL^1_1} \leq
\,&C\frac{e^t}{\nu(t)^{\frac{1}{2}}}\|f_0\|_{\rmL^1_1}+C\int_0^t
\frac{\rme^{t-s}}{\nu(t-s)^{\frac{1}{2}}}\|f(s)\|_{\infty} \|f(s)\|_{\rmL^1_1}\,\rmd s\\
&
+C\int_0^t\frac{\rme^{(t-s)}}{\nu(t-s)^{\frac{1}{2}}}\|f(w\cdot\nabla_w
f)(s)\|_{\rmL^1_1}\,\rmd s
\end{align*}
where by taking $p\geq 2$ and by interpolation as in
\eqref{eq:1norm_wf2}, we have
\begin{eqnarray*}
\|f(w\cdot\nabla_w f)\|_{\rmL^1_1}&\leq&
\||w|^{\frac{1}{2}}f\|_2 \||w|^{\frac{1}{2}}|\nabla_w f|\|_2\\
&\leq&
\|f\|_{\rmL_1^1}^{\frac{p-2}{2(p-1)}}\|f\|_{\rmL_m^p}^{\frac{p}{2(p-1)}}
\|\nabla_w f\|_{\rmL_1^1}^{\frac{p-2}{2(p-1)}}\|\nabla_w
f\|_{\rmL_m^p}^{\frac{p}{2(p-1)}}\\ &\leq&
\frac{\|f\|_{X_T}^2}{\nu(t)^{1/2}}.
\end{eqnarray*}
Putting together the last estimates, we deduce
\begin{equation}
\nu(t)^{\frac{1}{2}}\|\nabla_v f(t)\|_{\rmL^1_1}\leq
C_1^3(T,N,p)\|f_0\|_{\rmL^1_1}+C_2^3(T,N,p,r)\|f\|_{X_T}^2 \label{eq:der2}
\end{equation}
with $C_1^3$ and $C_2^3$ increasing functions of $T$,
for any $0< t\leq T$. From \eqref{eq:der1} and \eqref{eq:der2} and
all the estimates of the previous section, we finally get
$$
\|f\|_{X_T}\leq C_1(T,N,p)\|f_0\|_{\U}+C_2(T,N,p,r)\|f\|_{X_T}^2
$$
for any $T>0$. From these estimates and proceeding as at the end of
the previous section, it is easy to show that we have uniform
estimates in $X_T$ of the iteration sequence and the convergence
of the iteration sequence in the space $X_T$. From the uniqueness
obtained in the previous section, we conclude that the solution
obtained in this new procedure is the same as before and lies in
$X_T$. Summarizing, we have shown:

\begin{theorem}\label{exisder}
Let $m\geq 1$, $p>N$, $p\geq 2$, and $f_0\in\U$. Then there exists
$T>0$ depending only on the norm of the initial condition $f_0$ in
$\U$ such that \eqref{eq:fpf} has a unique solution in
$\calC([0,T];\U)$ with $f(0)=f_0$ and velocity gradients verifying that
$t\mapsto \nu(t)^{\frac{1}{2}} |\nabla_v f(t)| \in
BC((0,T),({\rmL_m^p}\cap{\rmL^1})(\R^N))$.
\end{theorem}

%
%
\subsection{Properties of the solutions}\label{s:aPS}

As \eqref{eq:fpf} belongs to the general class of convection-diffusion
equation, it enjoys several classical properties which we gather in
this section. The proof of these results uses classical approximation
arguments, see \cite{bib:ez,bib:v} for instance. Since these arguments are
somehow standard we will only give the detailed proof of the
$\rmL^1$-contraction property below.

\begin{lemma}[Positivity and Boundedness]\label{l:positivity}
Let $f\in X_T$ be the solution of the Cauchy problem
\eqref{eq:fpf} with initial condition $f_0\in\U$. If $0\leq f_0\leq 1$
in $\R^N$, then $0\leq f(t)\leq 1$ for any $0<t\leq T$.
\end{lemma}

\begin{lemma}[$L^1$-Contraction and Comparison Principle]\label{l:comparison}
Let $f\in X_T$ and $g\in X_T$ be the solutions of the Cauchy problem
\eqref{eq:fpf} with respective initial data $f_0\in\U$ and $g_0 \in \U$. Then
\begin{equation}
\label{spirou}
\|f(t)-g(t)\|_1\leq \|f_0-g_0\|_1
\end{equation}
for all $0<t\leq T$. Furthermore, if $f_0\le g_0$ then $f(t,v)\leq g(t,v)$
for all $0<t\leq T$ and $v\in\R^N$.
\end{lemma}
\proof
Since $f$ and $g$ solve (\ref{eq:fpf}),
\begin{equation}\label{kkk}
\frac{d}{d t}(f-g)=\Delta_v(f-g) + \nabla_v(v(f-g)) - \nabla_v(v(f^2-g^2))
\end{equation}
holds. We will obtain this result from the time evolution of
$|f-g|_{\e}$ where $|\cdot|_{\e}$ is the primitive vanishing at zero of $\sign_{\e}(s)$, the latter being an increasing smooth approximation of the
$\sign$ function defined by $\sign(s)=1$ if $s>0$, $\sign(0)=0$ and $\sign(s)=-1$ if $s<0$. Multiplying both sides of equation \eqref{kkk} by
$\zeta_n(v)\sign_{\e}(f-g)$ and integrating over $\R^N$,
where $\zeta_n\in\calC_0^{\infty}(\R^N)$ is a cut-off function
satisfying $0\leq\zeta_n\leq 1$, $\zeta_n(v)=1$ if $|v|\leq n$,
$\zeta_n(v)=0$ if $|v|\geq 2n$, and
$|\nabla_v\zeta_n|\leq\frac{1}{n}$, we obtain
\begin{align*}
\frac{d}{d t}\int_{\R^N}\zeta_n(v)|f-g|_{\e}\,\rmd v
\leq\,&-\int_{\R^N}\zeta_n(v)\sign_{\e}'(f-g)(v\cdot\nabla_v(f-g))(f-g)\,\rmd v\\
&
+\int_{\R^N}\zeta_n(v)\sign'_{\e}(f-g)(v\cdot\nabla_v(f-g))(f^2-g^2)\,\rmd
v\\ &
-\int_{\R^N}\nabla_v\zeta_n\sign_{\e}(f-g)(\nabla_v(f-g)+v(f-g-(f^2-g^2)))\,\rmd v\\
=&
-\int_{R^N}\zeta_n(v)(v\cdot\nabla_v((f-g)\sign_{\e}(f-g)-|f-g|_{\e}))\,\rmd
v\\ &
+\int_{\R^N}\zeta_n(v)(f+g)(v\cdot\nabla_v((f-g)\sign_{\e}(f-g)-|f-g|_{\e}))\,\rmd
v\\ &
-\int_{\R^N}\nabla_v\zeta_n\sign_{\e}(f-g)(\nabla_v(f-g)+v(f-g-(f^2-g^2)))\,\rmd
v.
\end{align*}
Integrating by parts, we finally get
\begin{align*}
\frac{d}{d t}\int_{\R^N}\zeta_n(v)|f-g|_{\e}\,\rmd v \leq\,&
\int_{\R^N}\div_v(v\zeta_n(v))((f-g)\sign_{\e}(f-g)-|f-g|_{\e})\,\rmd v\\
&
-\int_{\R^N}\div_v(\zeta_n(v)v(f+g))((f-g)\sign_{\e}(f-g)-|f-g|_{\e})\,\rmd
v\\& +\frac{1}{n}\int_{\R^N}|\nabla_v(f-g)+v(f-g-(f^2-g^2))|\,\rmd
v.
\end{align*}
For every $n$, the first two integrals become zero as
$\e\to 0$, since $f$ and $g$ are in $X_T$ whence $f(t),g(t)\in
L^1_1\cap L^\infty(\R^N)$ and $\nabla_v f(t),\nabla_v g(t)\in
L^1_1(\R^N)$ for any $0<t\leq T$, allowing for a Lebesgue
dominated convergence argument. We have that $\nabla_v
f+vf(1-f)\in L^1(\R^N)$ and $\nabla_v g+vg(1-g)\in L^1(\R^N)$ for
any $0<t\leq T$, and thus the third integral vanishes as
$n\to\infty$, getting finally
\begin{equation}
\frac{d}{d t}\int_{\R^N}|f-g|\,\rmd v\leq 0
\end{equation}
which concludes the proof of the first assertion of the lemma.
\endproof

Similar arguments show the conservation of mass.

\begin{lemma}[Mass Conservation]\label{l:conservation}
\label{l1} Let $f\in X_T$ be the solution of the Cauchy problem
\eqref{eq:fpf} with non-negative initial condition $f_0\in \U$, then the
$L^1$-norm of $f$ is conserved, i.e. $\|f(t)\|_1=\|f_0\|_1$ for
all $t\in [0,T]$.
\end{lemma}

Finally, we establish time dependent bounds on moments of
the solution to \eqref{eq:fpf}. More precisely, we will show that
moments increase at most as a polynomial on $t$. First, let us note that given
$a,b\geq 1$ and $f\in\rmL^1_{ab}(\R^N)\cap L^\infty(\R^N)$ then
\begin{equation}\label{eq:lpmbddbyl1pm}
\|f\|_{\rmL_a^b}\leq C\|f\|_{\rmL_{ab}^1}^{\frac{1}{b}}\|f\|_{\infty}^{1-\frac{1}{b}}.
\end{equation}
Indeed,
\begin{align*}
\|f\|_{\rmL_a^b}&=\left(\int_{\R^N}(1+|v|^a)^bf^b\rmd v\right)^{\frac{1}{b}} \leq\left(C\int_{\R^N}(1+|v|^{ab})f^b\rmd v\right)^{\frac{1}{b}}\nonumber\\
&\leq\left(C\|f\|_{\infty}^{b-1}\int_{\R^N}(1+|v|^{ab})f\rmd v\right)^{\frac{1}{b}}
=C\|f\|_{\rmL^1_{ab}}^{\frac{1}{b}}\|f\|_{\infty}^{1-\frac{1}{b}}.
\end{align*}
In particular, $(\rmL_{mp}^1\cap L^\infty)(\R^N)\subset \U$.

\noindent We next define $\lceil\g\rceil$ to be the smallest integer larger or equal than $\g$.

\begin{lemma}[Moments Bound]\label{l:momentgrowth}
Let $f\in X_T$ be the solution of the Cauchy problem \eqref{eq:fpf} with
initial condition $f_0\in \rmL_{mp}^1(\R^N)$ for some $m\ge 1$, $p>N$, $p\ge 2$, and satisfying $0\le f_0\le 1$. Then, for $0\leq t\leq T$
and $1\leq\g\leq mp/2$ the $2\g$-moment of $f(t)$ is bounded
by a polynomial $P_{\lceil\g\rceil}(t)$ of degree
$\lceil\g\rceil$, which depends only on the moments of $f_0$.
\end{lemma}

\proof We will prove it by induction on $\g$. First, we will
see that the second moment is bounded, and afterward that we can
bound every moment of order smaller than $pm$ in terms of a
$\g_*^{\text{th}}$ moment with $0<\g_*\leq 2$, which can in turn
be bounded in terms of the second moment.

Let $(\zeta_n)_{n\ge 1}$ be a sequence of smooth cut-off functions satisfying $0\le\zeta_n\le 1$, $\zeta_n(v)=1$ if $|v|\le n$, $\zeta_n(v)=0$ if $|v|\ge 2n$, $|\nabla_v\zeta_n|\leq 1/n$ and $|\D_v\zeta_n|\leq 1/n^2$. We multiply \eqref{eq:fpf} by $|v|^2\zeta_n(v)$ and integrate over $\R^N$ to get
\begin{align*}
\frac{d}{dt}\int_{\R^N}\zeta_n&(v)|v|^2f \rmd v=\int_{\R^N}\zeta_n(v)|v|^2 \D_v f \rmd v +\int_{\R^N}\zeta_n(v)|v|^2\div_v(vf(1-f)) \rmd v\\
&\leq \int_{\R^N}\left[\D_v\zeta_n|v|^2+4\nabla_v\zeta_n v+ 2N\zeta_n\right] f \rmd v + \int_{\R^N}|\nabla_v\zeta_n||v|^3f(1-f)\rmd v \\&\quad -2\int_{\R^N}\zeta_n|v|^2f\rmd v + 2\int_{\R^N}\zeta_n|v|^2f^2\rmd v\\
&\leq 5\int_{n<|v|<2n}f\rmd v + 2N\int_{\R^N}\zeta_n f\rmd v +\int_{n<|v|<2n}|v|^2f\rmd v .
\end{align*}
Now, letting $n\to\infty$ and noticing that $f\mathbbm{1}_{\{n<|v|<2n\}}$ and
$|v|^2f\mathbbm{1}_{\{n<|v|<2n\}}$ converge pointwise to zero and
are bounded by $f$ and $|v|^2f$ respectively with $f\in X_T$, we infer from the Lebesgue dominated convergence theorem that the first and the
last integrals converge to zero. Finally, integrating in time, we get
\begin{equation}
\int_{\R^N}|v|^2f(t,v)\rmd v \leq \int_{\R^N}|v|^2f_0(v)\rmd v +
2NMt
\end{equation}
for all $0\leq t\leq T$. Now, for the moment $2\g$ we can see in
the same way
\begin{align*}
\frac{d}{dt}\int_{\R^N}\zeta_n&(v)|v|^{2\g}f \rmd v=\int_{\R^N}\zeta_n(v)|v|^{2\g} \D_v f \rmd v +\int_{\R^N}\zeta_n(v)|v|^{2\g}\div_v(vf(1-f)) \rmd v\\
&\leq \int_{\R^N}\left[\D_v\zeta_n|v|^{2\g}+4\g\nabla_v\zeta_n|v|^{2(\g-1)}v+ 2\g(2(\g-1)+N)|v|^{2(\g-1)}\zeta_n\right] f \rmd v \\&\quad + \int_{\R^N}|\nabla_v\zeta_n||v|^{2\g+1}f(1-f)\rmd v  -2\g\int_{\R^N}\zeta_n|v|^{2\g}f\rmd v + 2\g\int_{\R^N}\zeta_n|v|^{2\g}f^2\rmd v\\
&\leq C\int_{n<|v|<2n}|v|^{2(\g-1)}f\rmd v + 2\g(2(\g-1)+N)\int_{\R^N}\zeta_n |v|^{2(\g-1)}f\rmd v \\&\quad+\int_{n<|v|<2n}|v|^{2\g}f\rmd v
\end{align*}
and we again let $n$ go to infinity. If $2\g\leq mp$, the
previous argument ensures that only the second integral remains,
and integrating in time, we conclude
\begin{equation}
\int_{\R^N}|v|^{2\g}f(t,v)\rmd v \leq
\int_{\R^N}|v|^{2\g}f_0(v)\rmd v +
2\g(2(\g-1)+N)\int_0^t\!\!\int_{\R^N} |v|^{2(\g-1)}f(s,v)\rmd v\,
\rmd s
\end{equation}
for all $0\leq t\leq T$. Whence, if we assume by induction that
the hypothesis of the lemma holds true for the $2(\g-1)$-moment,
\begin{equation}
\int_{\R^N}|v|^{2\g}f(v,t)\rmd v\leq
\int_{\R^N}|v|^{2\g}f_0(v)\rmd v + 2\g(2(\g-1)+N)\,\int_0^t
P_{\lceil\g-1\rceil}(s)\,\rmd s
\end{equation}
for all $0\leq t\leq T$, defining by induction the polynomial
$P_{\lceil\g\rceil}$.
\endproof

\subsection{Global existence}\label{s:gEiU}

Given an initial condition $f_0\in\rmL_{mp}^1(\R^N)$, $p>N$, $p\geq 2$, $m\geq 1$ such that $0\leq f_0\leq 1$, we have $f_0\in\U$ and we have shown in the previous subsections that there exists a unique local solution of (\ref{eq:fpf}) on an interval $[0,T)$. In fact, we can extend this solution to be global in time. If there exists $T_{max}<\infty$ such that the solution does not exist out of $(0,T_{max})$, then the $\U$-norm of it shall go to infinity as $t$ goes to $T_{max}$; as we will see, that situation cannot happen.

Due to Lemma \ref{l:positivity}, we have that
$0\leq f(t,v)\leq 1$ for any $0\leq t< T$ and any $v\in\R^N$, and
thus a bound for the $L^{\infty}$-norm of $f(t)$. Also, the
conservation of the mass in Lemma \ref{l:conservation} together
with the positivity in Lemma \ref{l:positivity} provide us with a
bound for the $\rmL^1$-norm. Finally, due to \eqref{eq:lpmbddbyl1pm} and Lemma \ref{l:momentgrowth} the $L^p_m$-norm is also bounded on any finite time interval.

\begin{theorem}[Global Existence]\label{main}
Let $f_0\in \rmL_{mp}^1(\R^N)$, $p>N$, $p\geq 2$, $m\geq 1$ be such
that $0\leq f_0\leq 1$. Then the Cauchy problem \eqref{eq:fpf}
with initial condition $f_0$ has a unique solution defined in
$[0,\infty)$ belonging to $X_T$ for all $T>0$. Also, we have $0\leq f(t,v)\leq 1$, for all $t\geq 0$ and $v\in\R^N$ and $\|f(t)\|_1=\|f_0\|_1=M$ for all $t\geq 0$.
\end{theorem}

\begin{remark}
Note that for any $K>0$ we can consider \eqref{eq:fpf} restricted to
the cylinder $C_K:=[0,\infty)\times\{|v|\leq K\}$. Then, due to the
fact that the solutions to \eqref{eq:fpf} we have constructed are in $\rmL^{\infty}$, we can show that
the solution is indeed $\calC^{\infty}(C_K)$ by applying regularity results in {\rm \cite{bib:lsu}} for
quasilinear parabolic equations.
\end{remark}

\begin{corollary}\label{c:radialsolution}
If $f_0\in\rmL^1_{mp}(\R^N)\cap\rmL^{\infty}(\R^N)$ is a radially symmetric and non-increasing function (that is, $f_0(v)=\varphi_0(|v|)$ for some non-increasing function $\varphi_0$), then so is $f(t)$ for all $t\geq 0$, that is, $f(t,v)=\varphi(t,|v|)$ and $r\mapsto\varphi(t,r)$ is non-increasing for all $t\geq 0$. In addition, $\varphi$ solves
\begin{equation}
\frac{\p\f}{\p t}=\frac{1}{r^{N-1}}\frac{\p}{\p
r}\left(r^{N-1}\frac{\p\f}{\p
r}+r^N\f(1-\f)\right) \quad\mbox{ with } \quad \frac{\p\f}{\p r}(t,0)=0 \label{eq:fpfRad}
\end{equation}
and $\f(0,r)=\f_0(r)$.
\end{corollary}
\proof
The uniqueness part of Theorem \ref{main} and the rotational invariance of \eqref{eq:fpf} imply that $f(t)$ is radially symmetric for all $t\ge 0$. The other properties are proved by classical arguments, the monotonicity of $r\mapsto \f(t,r)$ being a consequence of the comparison principle applied to the equation solved by $\p\f/\p r$.
\endproof

\section{Asymptotic Behaviour}\label{s:AB}
\setcounter{equation}{0}

Now that we have shown that under the appropriate assumptions
equation (\ref{eq:fpf}) have a unique solution which is global in
time, we are interested in how does this solution behave when the
time is large. For that we will define an appropriate entropy
functional for the solution and study its properties.


\subsection{Associated Entropy Functional}
In this section, we will show that the solutions constructed above
satisfy an additional dissipation property, the entropy decay. For $g\in\U$ such that $0\leq g\leq 1$, we define the functional
\begin{equation}\label{eq:functionalv2}
H(g):=S(g)+E(g)
\end{equation}
with the entropy given by
\begin{equation}
S(g):=\int_{\R^N} s(g(v))\,\rmd v
\end{equation}
where
\begin{equation}
s(r):=(1-r)\log(1-r)+r\log(r) \leq 0,\quad r\in[0,1],
\end{equation}
and the kinetic energy given by
\begin{equation}
E(g):=\frac{1}{2}\int_{\R^N}|v|^2g(v)\,\rmd v.
\end{equation}
We first check that $H(g)$ is indeed well defined and establish a control of the entropy in terms of the kinetic energy.

\begin{lemma}[Entropy Control]
For $\e\in(0,1)$, there exists a positive constant $C_{\e}$ such that
\begin{equation}\label{eq:relSE}
0\leq -S(g)\leq \e E(g)+ C_{\e}
\end{equation}
for every $g\in\rmL^1_2(\R^N)$ such that $0\leq g\leq 1$.
\end{lemma}
\proof
For $\e\in(0,1)$ and $v\in\R^N$, we put $z_{\e}(v):=1/(1+\rme^{\e|v|^2/2})$. The convexity of $s$ ensures that
\begin{align*}
s(g(v))-s(z_{\e}(v)) &\geq s'(z_{\e}(v))(g(v)-z_{\e}(v))\\
-s(z_{\e}(v))+s(g(v)) &\geq  \log\left(\frac{z_{\e}(v)}{1-z_{\e}(v)}\right)(g(v)-z_{\e}(v))
\end{align*}
for $v\in\R^N$. Since $z_{\e}(v)/(1-z_{\e}(v))=\rme^{-\e|v|^2/2}$, we end up with
\begin{align}
-s(g(v))&\leq\frac{\e|v|^2}{2}g(v)-s(z_{\e}(v))-\frac{\e|v|^2}{2}z_{\e}(v)\nonumber\\
&=\frac{\e|v|^2}{2}g(v)+(1-z_\e(v))\log\left(1+\rme^{-\e|v|^2/2}\right)+z_{\e}(v)\log\left(1+\rme^{-\e|v|^2/2}\right)
\nonumber\\
&\leq \frac{\e|v|^2}{2}g(v)+\rme^{-\e|v|^2/2}\label{domination}
\end{align}
for $v\in\R^N$, where we used $\log(1+a)\leq a$ for $a\geq 0$ and
$0\leq z_\e\leq 1$. Integrating the previous inequality yields
(\ref{eq:relSE}).
\endproof
\mbox{}\\
We next recall that $F_M$ is the unique Fermi-Dirac equilibrium state
satisfying $\|F_M\|_1=M:=\|f_0\|_1$; then we can introduce the next property
for $H$.

\begin{lemma}[Entropy Monotonicity]\label{l:hmonoton}
Assume that $f$ is the solution to the Cauchy problem
\eqref{eq:fpf} with initial condition $f_0$ in $\rmL_{mp}^1(\R^N)$ for some $p>max(N,2)$, $m\geq 1$ and satisfying $0\leq f_0\leq 1$. Then, the function $H$
is a non-increasing function of time satisfying for all $t>0$ that
\begin{equation}\label{eq:Hf0leqhftleqHF}
H(f_0) \geq H(f(t))\geq H(F_M) \quad\mbox{ with }\quad M:=\|f_0\|_1.
\end{equation}
\end{lemma}

\proof We first give a formal proof of the time monotonicity of $H(f)$ and supply additional details at the end of the proof. First of all, we observe that we can formulate
(\ref{eq:fpf}) as
$$
\frac{\p f}{\p
t}=\div_v\left[f(1-f)\nabla_v\left(s'(f)+\frac{|v|^2}{2}\right)\right].
$$
We multiply the previous equation by $s'(f)+|v|^2/2$ and integrate over $\R^N$ to obtain that
\begin{equation}\label{eq:ddtHleq0}
\frac{d}{dt}H(f)=-\int_{\R^N} f(1-f)|v+\nabla_vs'(f)|^2\,\rmd
v\leq 0.
\end{equation}

Consequently, the function $t\longrightarrow H(f(t))$ is a
non-increasing function of time, whence the first inequality in
(\ref{eq:Hf0leqhftleqHF}). To prove the second inequality, we
observe that the convexity of $s$ entails that
\begin{align*}
s(f(t,v))-s(F_M(v)) &\geq s'(F_M(v))(f(t,v)-F_M(v))\\
s(F_M(v))-s(f(t,v)) &\leq  \left(\log \b(M) +\frac{|v|^2}{2}\right)(f(t,v)-F_M(v))
\end{align*}
for $(t,v)\in [0,\infty)\times\R^N$. The second inequality in
(\ref{eq:Hf0leqhftleqHF}) now follows from the integration of the
previous inequality over $\R^N$ since
$\|F_M\|_1=\|f(t)\|_1$ by Lemma \ref{l:conservation}.

We shall point out that, in order to justify the previous
computations leading to the time monotonicity of the entropy, one should first start with an initial condition $f_0^\varepsilon$, $\varepsilon\in
(0,1)$, given by
$$
f_0^\varepsilon(v) = \max{\left\{\min{\left\{ f_0(v) , \frac{1}{1 +
\varepsilon
      e^{\vert v\vert^2/2}}\right\}} , \frac{\varepsilon}{\varepsilon
  + e^{\vert v\vert^2/2}} \right\}} \in \left[ \frac{\varepsilon}{\varepsilon +
    e^{\vert v\vert^2/2}} , \frac{1}{1 + \varepsilon e^{\vert
      v\vert^2/2}} \right]\,, \, v\in\R^N.
$$
Owing to the comparison principle (Lemma \ref{l:comparison}), the corresponding solution
$f^\varepsilon$ to (\ref{eq:fpf}) satisfies
$$
0 < \frac{\varepsilon}{\varepsilon + e^{\vert v\vert^2/2}} \le
f^\varepsilon(t,v) \le \frac{1}{1 + \varepsilon e^{\vert
v\vert^2/2}} < 1\,, \quad (t,v)\in (0,\infty)\times\R^N\,,
$$
for which the previous computations can be performed since the
solutions are immediately smooth and fast decaying at infinity for
all $t>0$, and thus $H(f^\e(t))\leq H(f^\e_0)$ for all $t\geq 0$.

Since $f_0^\varepsilon\to f_0$ in $\U$ and in $L^1_{mp}(\R^N)$ as
$\varepsilon\to 0$, it is not difficult to see that redoing all
estimates in subsections 2.1 and 2.2, we have continuous
dependence of solutions with respect to the initial data, and
thus, $f^\varepsilon$ converges towards $f$ in $X_T$ for any $T>0$.
Moreover, we have uniform bounds with respect to $\e$ of the moments in
finite time intervals using Lemma \ref{l:momentgrowth}. Direct
estimates easily show that $H(f^\e_0) \to H(f_0)$ as $\e\to 0$.

Let us now prove that $H(f^\e(t))\to H(f(t))$ as $\e\to 0$ for
$t>0$. Let us fix $R>0$.
Since $f^\varepsilon(t) \to f(t)$ in $L^1(\R^N)$ and we have
uniform estimates in $\e$ of moments of order $mp>2$ then
\begin{align*}
\Bigg{|}\int_{\R^N}|v|^2(f^\e(t) - f(t))\,\rmd v\Bigg{|} &\leq
\int_{|v|\geq R}|v|^2|f^\e(t) - f(t)|\,\rmd v
+\Bigg{|}\int_{|v|< R}|v|^2(f^\e(t) - f(t))\,\rmd v\Bigg{|}\\
&\leq \frac{1}{R^{mp-2}} \int_{|v|\geq R}|v|^{mp}(f^\e(t) +
f(t))\,\rmd v \\
&\qquad+ R^2 \left\|f^\e(t) - f(t)\right\|_1
\\
&\leq \frac{C(t)}{R^{mp-2}} + R^2 \left\|f^\e(t) - f(t)\right\|_1.
\end{align*}
Since the above inequality is valid for all $R>0$, we
conclude that $E(f^\e(t))\to E(f(t))$ as $\e\to0$. Now, taking
into account that $(1+|v|^2) f^\varepsilon(t) \to (1+|v|^2) f(t)$
in $L^1(\R^N)$, we deduce that there exists $h\in L^1(\R^N)$ such
that $||v|^2 f^\e(t)|\leq h$ and $f^\e(t)\to f(t)$ a.e. in $\R^N$,
for a subsequence that we denote with the same index. Using
inequality \eqref{domination}, we deduce that
$$
0\leq -s(f^\e(t,v))\leq\frac{1}{4}h(v)+\rme^{-|v|^2/4}\in
\rmL^1(\R^N)
$$
and that $-s(f^\e(t,v))\to -s(f(t,v))$ a.e. in $\R^N$. Thus, by
the Lebesgue dominated convergence theorem, we finally deduce that
$S(f^\e(t))\to S(f(t))$ as $\e\to0$. The convergence as $\e\to 0$ of $S(f^\e(t))$ to $S(f(t))$ is actually true for the whole family (and not only for a subsequence) thanks to the uniqueness of the limit. As a consequence, we showed
$H(f^\e(t))\to H(f(t))$ as $\e\to 0$ and passing to the limit
$\e\to0$ in the inequality $H(f^\e(t))\leq H(f^\e_0)$, we get the
desired result.
\endproof
\mbox{}

Now, it is easy to see the existence of a uniform in time bound
for the kinetic energy $E(f(t))$, or equivalently, of the
solutions in $\rmL^1_2(\R^N)$. If we take equations
(\ref{eq:functionalv2}), (\ref{eq:relSE}) (with $\e=1/2$) and
(\ref{eq:Hf0leqhftleqHF}) we get that
$$
E(f(t))=H(f(t))-S(f(t))\leq \frac{1}{2}E(f(t))+C_{1/2}+H(f_0)
$$
for $t\geq 0$ whence
\begin{equation}\label{eq:bnd2ndmoment}
E(f(t))\leq 2\big(C_{1/2}+H(f_0)\big).
\end{equation}

\subsection{Convergence to the Steady State}
\begin{theorem}[Convergence]\label{p:convergence}
Let $f$ be the solution to the Cauchy problem
\eqref{eq:fpf} with initial condition $f_0$ in $\rmL_{mp}^1(\R^N)$,
$p>\max(N,2)$, $m\geq 1$ satisfying $0\leq f_0 \leq 1$. Then
$\{f(t)\}_{t\geq 0}$ converges strongly in $\rmL^1(\R^N)$ towards
$F_M$ as $t\to\infty$ with $M:=\|f_0\|_1$.
\end{theorem}

For the proof, we first need a technical lemma.

\begin{lemma}\label{tech2} Let $f$ be the solution to the Cauchy problem
\eqref{eq:fpf} with initial condition $f_0$ in $\rmL_{mp}^1(\R^N)$,
$p>\max(N,2)$, $m\geq 1$ satisfying $0\leq f_0 \leq 1$. If $A$ is a measurable subset of $\R^N$, we have
\begin{equation}\label{z3}
\int_0^\infty \left( \int_A \left\vert v f(1-f) + \nabla_v f
  \right\vert dv \right)^2 dt \le H(F_M) \sup_{t\geq 0} \left\{ \int_A
f(t,v) \rmd v \right\}
\end{equation}
\end{lemma}
\proof Owing to the second inequality in (\ref{eq:Hf0leqhftleqHF})
and the finiteness of $H(f_0)$, we also infer from
(\ref{eq:ddtHleq0}) that $(t,v)\longmapsto f (1-f) \left| v +
\nabla_v s'(f) \right|^2$ belongs to $L^1((0,\infty)\times\R^N)$.
Working again with the regularized solutions $f^\varepsilon$, it
then follows from Lemma \ref{l:conservation} and the
Cauchy-Schwarz inequality that, if $A$ is a measurable subset of
$\R^N$, we can compute
\begin{align*}
\int_0^\infty \Bigg{(} \int_A \vert v
f^\varepsilon(1-f^\varepsilon) + &\nabla_v f^\varepsilon
  \vert \rmd v\Bigg{)}^2 \rmd t \\&=  \int_0^\infty \left( \int_A
  \frac{\left\vert v f^\varepsilon(1-f^\varepsilon) + \nabla_v f^\varepsilon \right\vert}
{\left( f^\varepsilon(1-f^\varepsilon)
    \right)^{1/2}} \left( f^\varepsilon (1-f^\varepsilon) \right)^{1/2} \rmd v \right)^2 \rmd t \\
& \le  \int_0^\infty \left( \int_A \frac{\left\vert v
f^\varepsilon(1-f^\varepsilon) +
      \nabla_v f^\varepsilon \right\vert^2}{f^\varepsilon (1-f^\varepsilon)} \rmd v \right) \left( \int_A f^\varepsilon
  (1-f^\varepsilon) \rmd v \right) \rmd t ,
\end{align*}
and thus,
\begin{align*}
\int_0^\infty \Bigg{(} \int_A \vert v
f^\varepsilon(1-f^\varepsilon) + &\nabla_v f^\varepsilon
  \vert \rmd v\Bigg{)}^2 \rmd t \\
& \le  \sup_{t\ge 0} \left\{ \int_A f^\varepsilon(t,v) \rmd v
\right\} \int_0^\infty \int_A
f^\varepsilon (1-f^\varepsilon) \left[ v + \nabla_v s'(f^\varepsilon) \right]^2 \rmd vdt \\
& \le  H(F_{M^\e}) \sup_{t\ge 0} \left\{ \int_A
f^\varepsilon(t,v) \rmd v \right\} \,.
\end{align*}
Here, $M^\e:=\|f_0^\e\|_1$ so that $F_{M^\e}$ is the Fermi-Dirac distribution with the mass of the regularized initial condition $f_0^\e$. It is easy to check that
$H(F_{M^\e})\to H(F_M)$ as $\e\to0$ since $M^\e\to M$ as $\e\to0$.
Passing to the limit as $\varepsilon\to 0$, $f^\varepsilon\to f$ in
$X_T$ for any $T>0$, and thus we get the conclusion.
\endproof
\

\par{\it \noindent Proof of Theorem \ref{p:convergence}.-}
We first establish that
\begin{equation}
\label{z6}
\{f(t)\}_{t\ge 0} \;\;\mbox{ is bounded in }\;\; \rmL_2^1(\R^N)\cap
\rmL^{\infty}(\R^N)\,.
\end{equation}
{F}rom \eqref{eq:bnd2ndmoment} and Theorem \ref{main}, it is straightforward that $E(f(t))$ is bounded in $[0,\infty)$. Recalling the mass conservation, the boundedness of $\{f(t)\}_{t\ge 0}$ in $\rmL_2^1(\R^N)\cap \rmL^{\infty}(\R^N)$ follows.

We next turn to the strong compactness of $\{f(t)\}_{t\ge 0}$ in
$\rmL^1(\R^N)$. For that purpose, we put $R(t,v) := v f(t,v) (1-f(t,v))$
for $(t,v)\in (0,\infty)\times\R^N$ and deduce from Theorem \ref{main} and (\ref{z6}) that
\begin{equation}
\label{z7} \sup_{t\ge 0} \left( \Vert R(t)\Vert_1 + \Vert
R(t)\Vert_2^2 \right) \le 2\ \sup_{t\ge 0} \int_{\R^N} (1+|v|^2)
f(t,v) \rmd v < \infty\,.
\end{equation}
Denoting the linear heat semigroup on $\R^N$ by $(e^{t\Delta})_{t\ge
  0}$, it follows from (\ref{eq:fpf}) that $f$ is given by the Duhamel
formula
\begin{equation}
\label{z8}
f(t) = e^{t\Delta} f_0 + \int_0^t \nabla_v e^{(t-s)\Delta} R(s)
\rmd s\,,\quad t\ge 0\,.
\end{equation}
It is straightforward to check by direct Fourier transform techniques that
\begin{eqnarray*}
\Vert e^{t\Delta} g\Vert_{\dot{H}^\alpha} \le C(\alpha)\ \min{\left\{
    t^{-\alpha/2} \Vert g\Vert_2 , t^{-(2\alpha+N)/4} \Vert
    g\Vert_1 \right\}}
\end{eqnarray*}
for $t\in (0,\infty)$, $g\in \rmL^1(\R^N)\cap \rmL^2(\R^N)$ and $\alpha\in
[0,2]$ with
$$
\Vert g\Vert_{\dot{H}^\alpha} := \left( \int_{\R^N}
|\xi|^{2\alpha} \left|
    \widehat{g}(\xi)\right|^2 d\xi\right)^{1/2}
$$
and $\widehat{g}$ being the Fourier transform of $g$. Thus, we deduce from
(\ref{z8}) that, if $t\ge 1$ and $\alpha\in \left( (1-(N/2))^+,1
\right)$, we have

\begin{eqnarray*}
\Vert f(t)\Vert_{\dot{H}^\alpha} & \le & C(\alpha)
t^{-(2\alpha+N)/4} \Vert f_0\Vert_1 + C(\alpha+1)
\int_0^{t-1}
(t-s)^{-(2+2\alpha+N)/4} \Vert R(s)\Vert_1 \rmd s \\
& & +\ C(\alpha+1) \int_{t-1}^t (t-s)^{-(1+\alpha)/2} \Vert
R(s)\Vert_2 \rmd s \\
& \le & C \left( 1 + \int_1^t s^{-(2+2\alpha+N)/4} \rmd s +
\int_0^1
  s^{-(1+\alpha)/2} \rmd s \right)\\
& \le & C\,,
\end{eqnarray*}
thanks to the choice of $\alpha$. Consequently, $\{f(t)\}_{t\ge 1}$ is
also bounded in $\dot{H}^\alpha$ for $\alpha\in \left( (1-(N/2))^+,1
\right)$. Owing to the compactness of the embedding of $(\dot{H}^\alpha\cap \rmL_2^1)(\R^N)$ in $\rmL^1(\R^N)$, we finally conclude
that
\begin{equation}
\label{z9}
\{f(t)\}_{t\ge 0} \;\;\mbox{ is relatively compact in }\;\; \rmL^1(\R^N)\,.
\end{equation}

Consider now a sequence $\{t_n\}_{n\in\N}$ of positive real
numbers such that $t_n\to\infty$ as $n\to\infty$. Owing to
(\ref{z9}), there are a subsequence of $\{t_n\}$ (not relabelled)
and $g_\infty\in \rmL^1(\R^N)$ such that $\{f(t_n)\}_{n\in\N}$
converges towards $g_\infty$ in $L^1(\R^N)$ as $n\to\infty$.
Putting $f_n(t) = f(t_n+t)$, $t\in [0,1]$ and denoting by $g$ the
unique solution to (\ref{eq:fpf}) with initial datum $g_\infty$,
we infer from the contraction property \eqref{spirou} that
\begin{equation}
\label{z10}
\lim_{n\to\infty} \sup_{t\in [0,1]} \Vert f_n(t) - g(t)\Vert_1 = 0\,.
\end{equation}
Next, on one hand, we deduce from the proof of Lemma \ref{tech2}
with $A=\R^N$ that $(t,v) \longmapsto v f(t,v) (1-f(t,v)) +
\nabla_v f(t,v)$ belongs to $\rmL^2((0,\infty);\rmL^1(\R^N))$. Since
$$
\int_0^1 \left( \int_{\R^N} \left\vert v f_n (1-f_n) + \nabla_v
f_n
  \right\vert \rmd v \right)^2 \rmd t = \int_{t_n}^{t_n+1} \left( \int_{\R^N}
  \left\vert v f (1-f) + \nabla_v f \right\vert \rmd v \right)^2 \rmd t\,,
$$
we end up with
\begin{equation}
\label{z11} \lim_{n\to\infty} \int_0^1 \left( \int_{\R^N}
\left\vert v f_n (1-f_n) +
    \nabla_v f_n \right\vert \rmd v \right)^2 \rmd t = 0\,.
\end{equation}
On the other hand, it follows from the mass conservation and (\ref{z3}) that, if
    $A$ is a measurable subset of $\R^N$ with finite measure $|A|$,
    we have
$$
\int_0^1 \left( \int_A \left\vert v f_n (1-f_n) + \nabla_v f_n
  \right\vert \rmd v \right)^2 \rmd t \le H(F_M) |A|\,,
$$
which implies that $\{v f_n (1-f_n) + \nabla_v f_n\}_{n\in\N}$ is
weakly relatively compact in $\rmL^1((0,1)\times\R^N)$ by the
Dunford-Pettis theorem. Since $\{v f_n (1-f_n)\}_{n\in\N}$
converges strongly towards $v g (1-g)$ in $\rmL^1((0,1)\times\R^N)$
by (\ref{z6}) and (\ref{z10}), we conclude that $\{\nabla_v
f_n\}_{n\ge 0}$ is weakly relatively compact in
$L^1((0,1)\times\R^N)$. Upon extracting a further subsequence, we
may thus assume that $\{\nabla_v f_n\}_{n\ge 0}$ converges weakly
towards $\nabla_v g$ in $\rmL^1((0,1)\times\R^N)$. Consequently,
$$
\int_0^1 \int_{\R^N} \left\vert v g (1-g) + \nabla_v g \right\vert
\rmd v\,\rmd t \le \liminf_{n\to\infty} \int_0^1 \int_{\R^N}
\left\vert v f_n (1-f_n) + \nabla_v
  f_n \right\vert \rmd v\,\rmd t = 0
$$
by (\ref{z11}), from which we readily deduce that $v g (1-g) +
\nabla_v g = 0$ a.e. in $(0,1)\times\R^N$. Since $\Vert
g(t)\Vert_1=M$ for each $t\in [0,1]$ by Lemma
\ref{l:conservation} and (\ref{z10}), standard arguments allow us
to conclude that $g(t)=F_M$ for each $t\in [0,1]$. We have thus
proved that $F_M$ is the only possible cluster point in $L^1(\R^N)$
of $\{f(t)\}_{t\ge 0}$ as $t\to\infty$, which, together with the
relative compactness of $\{f(t)\}_{t\ge 0}$ in $\rmL^1(\R^N)$,
implies the assertion of Theorem~\ref{p:convergence}.
\endproof

By now, we have seen that the solution of
(\ref{eq:fpf}) with initial condition $f_0$ converges to the Fermi-Dirac distribution $F_M$ with the same mass as $f_0$ as
$t\rightarrow\infty$, but we are also interested in how fast this
happens. We will answer that question with the next result, which
was already proved in \cite{bib:crs} in the one dimensional case,
and easily extends to any dimension based on the existence and
entropy decay results established above.

\begin{theorem}[Entropy Decay Rate]\label{thm:EntropyDecayRate}
Let $f$ be the solution to the Cauchy problem
\eqref{eq:fpf} with initial condition $f_0$ in $\rmL_{mp}^1(\R^N)$,
$p>\max(N,2)$, $m\geq 1$ satisfying $0\leq f_0 \leq F_{M^*} \leq 1$
for some $M^*$. Then
\begin{equation}\label{ll1}
H(f(t))-H(F_M) \leq (H(f_0)-H(F_M))\mathrm{e}^{-2Ct}
\end{equation}
and
\begin{equation}\label{ll2}
\|f(t)-F_M\|_1 \leq C_2 (H(f_0)-H(F_M))^{1/2} \mathrm{e}^{-Ct}
\end{equation}
for all $t\geq 0$, where $C$ depends on $M^*$ and $M:=\|f_0\|_1$.
\end{theorem}

\proof Since $0\leq f_0\leq F_{M^*}$, then the initial condition
satisfies all the hypotheses of Theorems~\ref{main} and~\ref{p:convergence}. In order to show the exponential
convergence, we use the same arguments as in \cite{bib:crs}. We
first remark that the entropy functional $H$ coincides with the
one introduced in \cite{bib:cjmtu} for the nonlinear diffusion
equation
\begin{equation}
\frac{\p g}{\p t}=\div_x\left[g \nabla_x\left( x + h(g)
\right)\right] \label{eq:nonlindif}
\end{equation}
for the function $0\leq g(t,x)\leq 1$, $x\in\R$, $t>0$, where
$h(g) =  s'(g) = \log g -\log (1-g)$. Let us point out that the
relation between the entropy dissipation for the solutions of the
nonlinear diffusion equation \eqref{eq:nonlindif}, given by
$$
-D_0(g) = \frac{d}{dt}H(g) = -\int_{\R^N}g\left|x+\frac{\p}{\p
x} h(g)\right|^2\,dx ,
$$
and the entropy dissipation for the solutions of \eqref{eq:fpf},
given by \eqref{eq:ddtHleq0}, is the basic idea of the proof.
Indeed, one can check that, once restricted to the range $f\in (0,1)$, $h(f)$ verifies the hypotheses of the Generalized Logarithmic Sobolev
Inequality \cite[Theorem~17]{bib:cjmtu}. The Generalized Logarithmic Sobolev
Inequality then asserts that
\begin{equation}\label{eq:cjmtu01}
H(g)-H(F_M) \leq \frac{1}{2}D_0(g)
\end{equation}
for all integrable positive $g$ with mass $M$ for which the
right-hand side is well-defined and finite. We can now, by the
same regularization argument as before, compare the entropy
dissipation $D(f)=-\frac{d}{d t}H(f)$ of equation \eqref{eq:fpf}
and the one $D_0(f)$ of equation \eqref{eq:nonlindif}. Thanks to
Lemma \ref{l:comparison} we have $f(t,v)\leq F_{M^*}(v)\leq
(\beta(M^*)+1)^{-1}$ a.e. in $\R^N$, and thus
\begin{equation}
D(f)=\int_{\R^N}f(1-f)\left|v+\nabla_vh(f)\right|^2 dv \geq
C\int_{\R^N}f \left|v+\nabla_vh(f)\right|^2 dv
\end{equation}
where $C=1-(\beta(M^*)+1)^{-1}$. Applying the Generalized Logarithmic Sobolev Inequality (\ref{eq:cjmtu01}) to the solution $f$ and taking
into account the previous estimates, we conclude
\begin{equation}
H(f(t)) - H(F_M) \leq (2C)^{-1} D(f(t)).
\end{equation}
Finally, coming back to the entropy evolution:
$$
\frac{d}{dt}\left[H(f(t))-H(F_M)\right]=-D(f(t)) \leq
-2C\left[H(f(t))-H(F_M)\right],
$$
and the result follows from Gronwall's lemma. The convergence in
$\rmL^1$ is obtained by a Csisz\'ar-Kullback type inequality proven in
\cite[Corollary 4.3]{bib:crs}, its proof being valid for any space
dimension. It is actually a consequence of a direct application of
the Taylor theorem to the relative entropy $H(f)-H(F_M)$ giving:
$$
\|f-F_M\|_1^2 \leq 2 M (H(f)-H(F_M)) .
$$
\endproof

\subsection{Propagation of Moments and Consequences}\label{s:PoM}

There is a large gap between Theorem~\ref{p:convergence} which only provides the $L^1$-convergence to the equilibrium and Theorem~\ref{thm:EntropyDecayRate} which warrants an exponential decay to zero of the relative entropy for a restrictive class of initial data. This last section is devoted to an intermediate result where we prove the convergence to zero of the relative entropy but without a rate for a larger class of initial data than in Theorem~\ref{thm:EntropyDecayRate}.

\begin{lemma}[Time independent bound for Moments]\label{l:bddmoments}
Let $g_0\in\rmL_{mp}^1(\R^N)$ with $m\geq 1$, $p>\max{(p,2)}$ such that $0\leq g_0\leq 1$, and assume further that $g_0$ is a radially symmetric and
non-increasing function, i.e., there is a non-increasing function $\f_0$ such that $g_0(v)=\f_0(|v|)$ for $v\in\R^N$. Then, for the unique solution $g$ of the Cauchy problem \eqref{eq:fpf} with initial condition $g_0$, the control of moments propagates in time, i.e., there exists $C>0$ depending on $N$ and
$g_0$, but not on time, such that
\begin{equation}\label{eq:cuamoment}
\lim_{R\to\infty}\sup_{t\geq 0}\int_{\{|v|\geq R\}} |v|^{mp} g(t,v)\rmd v=0.
\end{equation}
\end{lemma}

\proof We have already seen in Corollary \ref{c:radialsolution} the existence and uniqueness of $g$ and that $g(t,v)=\f(t,|v|)$ for $t\ge 0$ and $v\in\R^N$ for some function $\f$ such that $r\mapsto\f(t,r)$ is non-increasing. Furthermore, we have that its moments are given by
\begin{equation}\label{eq:radeqMass}
M:=\int_{\R^N} g(t,v)\,\rmd
v=N\omega_N\int_0^{\infty}r^{N-1}\f(t,r)\,\rmd r
\end{equation}
and
\begin{equation}\label{eq:radeqMoment}
\int_{\R^N} |v|^{mp} g(t,v)\,\rmd
v=N\omega_N\int_0^{\infty}r^{N+mp-1}\f(t,r)\,\rmd r
\end{equation}
for $t\geq 0$, where $\omega_N$ denotes the volume of the unit ball of
$\R^N$.

Next, since $|v|^{mp} g_0\in\rmL^1(\R^N)$, the map
$v\mapsto |v|^{mp}$ belongs to $\rmL^1(\R^N;g_0(v)\,dv)$ and a
refined version of de la Vall\'ee-Poussin theorem
\cite{bib:dm,bib:ch} ensures that there is a non-decreasing,
non-negative and convex function
$\psi\in\calC^{\infty}([0,\infty))$ such that $\psi(0)=0$, $\psi'$
is concave,
\begin{equation}\label{eq:Philimandf}
\lim_{r\rightarrow\infty}\frac{\psi(r)}{r}=\infty\quad\mbox{and}
\quad \int_{\R^N}\psi(|v|^{mp})g_0(v)\,\rmd v<\infty .
\end{equation}
Observe that, since $\psi(0)=0$ and $\psi'(0)\geq 0$, the convexity of
$\psi$ and the concavity of $\psi'$ ensure that for $r\geq 0$ 
\begin{equation}\label{eq:relphi}
r\psi''(r)\leq\psi'(r)\qquad \mbox{and} \qquad\psi(r)\leq
r\psi'(r).
\end{equation}

Then, after integration by parts, it follows from (\ref{eq:fpfRad}) that
\begin{align}\label{eq:fpfRadwithPhi}
\frac{1}{mp}\frac{d}{d t}\int_0^{\infty}\psi(r^{mp})r^{N-1}\f\,\rmd
r&=-\int_0^{\infty}r^{mp-1}\psi'(r^{mp})\left( r^{N-1}\frac{\p\f}{\p
r} + r^N\f(1-\f)\right)\,\rmd r\nonumber\\&= I_1+I_2 ,
\end{align}
where
\begin{eqnarray*}
I_1&=&\int_0^{\infty}\f\left[(mp+N-2)r^{mp+N-3}\psi'(r^{mp})+mpr^{2mp+N-3}\psi''(r^{mp})\right]\rmd r\\
I_2&=&-\int_0^{\infty}r^{N+mp-1}\psi'(r^{mp})\f(1-\f)\,\rmd r.
\end{eqnarray*}

We now fix $R>0$ such that $\omega_N R^N\geq 4M$ and $R^2\geq
4(2mp+N-2)$, and note that due to the monotonicity of $\f$ with
respect to $r$ and (\ref{eq:radeqMass})-(\ref{eq:radeqMoment}) the
inequality
\begin{equation}\label{tech1}
M\geq N\omega_N\int_0^Rr^{N-1}\f\rmd r\geq \omega_N R^N\f(R)
\end{equation}
holds. Therefore, we first use the monotonicity of $\psi'$ and
$\f$ together with \eqref{tech1} to obtain
\begin{eqnarray*}
I_2&\leq&-\int_R^{\infty}r^{N+mp-1}\psi'(r^{mp})\f(1-\f)\,\rmd r\leq (\f(R)-1)\int_R^{\infty}r^{N+mp-1}\psi'(r^{mp})\f\,\rmd r\\
&\leq&\left(\frac{M}{\omega_NR^N}-1\right)\int_R^{\infty}r^{N+mp-1}\psi'(r^{mp})\f\,\rmd r
\leq-\frac{3}{4}\int_R^{\infty}r^{N+mp-1}\psi'(r^{mp})\f\,\rmd r\\
&\leq&\frac{3}{4}\int_0^Rr^{N+mp-1}\psi'(r^{mp})\f\,\rmd r-\frac{3}{4}\int_0^{\infty}r^{N+mp-1}\psi'(r^{mp})\f\,\rmd r\\
&\leq&\frac{3MR^{mp}\psi'(R^{mp})}{4N\omega_N}-\frac{3}{4}\int_0^{\infty}r^{N+mp-1}\psi'(r^{mp})\f\,\rmd
r.
\end{eqnarray*}

On the other hand, from
(\ref{eq:radeqMass}),(\ref{eq:radeqMoment}), (\ref{eq:relphi}),
\eqref{tech1} and the monotonicity of $\psi'$
\begin{eqnarray*}
I_1&\leq&(N+2m-2)\int_0^{\infty}r^{N+mp-3}\psi'(r^{mp})\f\,\rmd r\\
&\leq&(N+2mp-2)\psi'(R^{mp})R^{mp-2}\int_0^Rr^{N-1}\f\,\rmd r\\& & +\frac{N+2mp-2}{R^2}\int_R^{\infty}r^{N+mp-1}\psi'(r^{mp})\f\,\rmd r\\
&\leq&
\frac{(N+2mp-2)\psi'(R^{mp})R^{mp-2}M}{N\omega_N}+\frac{1}{4}\int_R^{\infty}r^{N+mp-1}\psi'(r^{mp})\f\,\rmd
r.
\end{eqnarray*}

Inserting these bounds for $I_1$ and $I_2$ in (\ref{eq:fpfRadwithPhi}) and using (\ref{eq:relphi}) we end up with
\begin{align*}
\frac{1}{mp}\frac{d}{d t}\int_0^{\infty}&\psi(r^{mp})r^{N-1}\f\,\rmd r\\
&\leq \frac{\psi'(R^{mp})MR^{mp-2}}{N\omega_N}\left(\frac{3R^2}{4}+N+2mp-2\right)-\frac{1}{2}\int_0^{\infty}r^{N+mp-1}\psi'(r^{mp})\f\,\rmd r\\
&\leq
\frac{\psi'(R^{mp})MR^{mp-2}}{N\omega_N}\left(\frac{3R^2}{4}+N+2mp-2\right)
- \frac{1}{2}\int_0^{\infty}r^{N-1}\psi(r^{mp})\f\,\rmd r.
\end{align*}

We then use the Gronwall lemma to
conclude that there exists $C>0$ depending on $N$, $M$, $m$, $p$, $g_0$ and
$\psi$ such that
\begin{equation*}
\sup_{t\geq 0}\int \psi(|v|^{mp})g(t,v)\rmd v\leq C
\end{equation*}
from which (\ref{eq:cuamoment}) readily follows by
(\ref{eq:Philimandf}).
\endproof

\begin{theorem}[Entropy Convergence]
Let $f$ be the solution of the Cauchy problem \eqref{eq:fpf} with
initial condition $f_0\in \rmL^1_{mp}(\R^N)$ such that there exists a
radially symmetric and non-increasing function $g_0 \in \rmL^1_{mp}(\R^N)$ with $0\leq f_0\leq g_0\leq 1$. Then $H(f)\to H(F_M)$ as $t\to\infty$ where $M=\|f_0\|_1$.
\end{theorem}

\proof
Due to \cite[Theorem 3]{bib:lw} we know that
\begin{align*}
|H(f(t))-H(F_M)|&\leq C\int_{\R^N}|v|^2|f(t,v)-F(v)|\rmd w\\
&\leq R^2\|f(t)-F\|_1+\sup_{t\geq 0}\int_{|v|\geq R}|v|^2|f(t)-F|\rmd v
\end{align*}
Now, Theorem \ref{thm:EntropyDecayRate} and Lemma \ref{l:bddmoments}
imply that $H(f(t))\to H(F_M)$ as $t\to\infty$.
\endproof

%
%
%
%
\appendix
\section[Appendix]{$\rmL_m^p$-bounds for the Fokker-Planck Operator}\label{ap:bounds}
\setcounter{equation}{0}
\def\theequation{A.\arabic{equation}}

Here we follow similar arguments as in \cite{bib:gw} to show some
bounds for $\|\p_{\alpha}\calF f (t)\|_{\rmL^p_m} $ which were
useful in the fixed point argument in Section~\ref{s:lEiU}. We recall
the well-known Young inequality: Let 
$g_1\in\rmL^r(\R^N)$, $g_2\in\rmL^q(\R^N)$ with $1\leq
p,r,q\leq\infty$ and $\frac{1}{p}+1=\frac{1}{r}+\frac{1}{q}$, then
$$
g_1\ast g_2\in\rmL^p(\R^N)\quad\text{and}\quad \|g_1\ast
g_2\|_p\leq\|g_1\|_r\,\|g_2\|_q.
$$

\begin{proposition}\label{prop:fpnucpqbounds}
Let $1\leq q\leq p\leq \infty$, $m\geq 0$ and $\alpha\in \N^N$. Then
for $t>0$,
\begin{equation}
\|\p_{\a}\calF(t)[f]\|_{\rmL^p_m}\leq
\frac{C\rme^{\left(\frac{N}{p'}+|\a|\right)t}}{\nu(t)^{\frac{N}{2}\left(\frac{1}{q}-\frac{1}{p}\right)
+\frac{|\alpha|}{2}}}\|f\|_{\rmL^q_m} .
\end{equation}
\end{proposition}

\proof
For all $\alpha\in \N^N$, we have
\begin{align}
\p_{\a}\calF(t,v)[f] &= \p^{\a}\int_{\R^N} \left(
\frac{\rme^{t N}}
{\left(2\pi\left(\rme^{2t}-1\right)\right)^{\frac{N}{2}}}
\rme^{-\frac{|\rme^{t}v-w|^2}
{2\left(\rme^{2t}-1\right)}}\right)f(w)\ \rmd w\nonumber\\
&=\p^{\a}\int_{\R^N}\left(
\frac{\rme^{2Nt}}{\left(2\pi\left(\rme^{2t}-1\right)\right)^{\frac{N}{2}}}
\rme^{-\frac{|\rme^{t}(v-w)|^2}{2\left(\rme^{2t}-1\right)}}
\right)f(\rme^{t}w)\ \rmd w\nonumber\\
&=\frac{\rme^{t \left(2N+|\a|\right)}}{\nu(t)^{\frac{N+|\a|}{2}}}
\int_{\R^N} \phi_{\a}\left(\frac{v-w}{\rme^{-t}\nu(t)^{1/2}}\right)f(\rme^{t}w)\ \rmd w
\end{align}
where
$$
\phi_{\a}(\chi)=\p^{\a}_{\chi}\left(\phi_0\right)(\chi)=\calP_{|\a|}(\chi)\phi_0(\chi)\text{,}
$$
being $\calP_{|\a|}(\chi)$ a polynomial of degree $|\a|$ which we
can recursively reckon by
$$
\calP_0(\chi)=1\text{, }
\calP_{|\a|}(\chi)=\calP_{|\a|-1}'(\chi)-\chi\calP_{|\a|-1}(\chi)\text{
and }\phi_0(\chi)=(2\pi)^{-\frac{N}{2}}\rme^{-\frac{|\chi|^2}{2}}.
$$

Since $1+|v|^m\leq C(1+|v-w|^m)(1+|w|^m)$, we deduce
\begin{align}
(1+&|v|^m)|(\p_{\a}\calF\ast f)(t)|\leq\nonumber\\&\leq
C\frac{\rme^{t \left(2N+|\a|\right)}}{\nu(t)^{\frac{N+|\a|}{2}}}
\int_{\R^N}
(1+|v-w|^m)\Bigg{|}\phi_{\a}\left(\frac{v-w}{\rme^{-t}\nu(t)^{1/2}}\right)\Bigg{|}(1+|w|^m)\Big{|}f(\rme^{t}w)\Big{|}\rmd
w\label{eq:dNucFPf}\text{.}
\end{align}
Then, we can write
\begin{align*}
\int_{\R^N}&(1+|v-w|^m)^r\Bigg{|}\phi_{\a}\left(\frac{v-w}{\rme^{-t}\nu(t)^{1/2}}\right)\Bigg{|}^r\rmd w
=C(I+II)
\end{align*}
with
\begin{align*}
I=\int\calP_{|\a|}^r\left(\frac{v-w}{\rme^{-t}\nu(t)^{1/2}}\right)\phi_0\left(\frac{v-w}{\rme^{-t}\nu(t)^{1/2}}\right)^r\rmd
w =
\frac{\nu(t)^{N/2}}{\rme^{Nt}}\int\calP_{|\a|}^r(\chi)\phi_0(\chi)^r=C_1\frac{\nu(t)^{N/2}}{\rme^{Nt}}
\end{align*}
and
\begin{align*}
II&=\int|v-w|^{mr}\calP_{|\a|}^r\left(\frac{v-w}{\rme^{-t}\nu(t)^{1/2}}\right)\phi_0\left(\frac{v-w}{\rme^{-t}\nu(t)^{1/2}}\right)^r\rmd
w \\&=
\frac{\nu(t)^{(N+mr)/2}}{\rme^{(N+mr)t}}\int|\chi|^{mr}\calP_{|\a|}^r(\chi)\phi_0(\chi)^r=C_2\frac{\nu(t)^{(N+mr)/2}}{\rme^{(N+mr)t}}.
\end{align*}
whence
\begin{equation}\label{eq:phiweightlpnorm}
\frac{\rme^{Nt}}{\nu(t)^{N/2}}\int_{\R^N}(1+|v-w|^m)^r\Bigg{|}\phi_{\a}
\left(\frac{v-w}{\rme^{-t}\nu(t)^{1/2}}\right)\Bigg{|}^r\rmd w
\leq C .
\end{equation}

On the other hand, we get
\begin{align}
\Big{\|}(1+|w|^m)\bigg{|}f(\rme^{t}w)\bigg{|}\Big{\|}_p
&=\left(\int (1+|w|^m)^p\Big{|}f(\rme^{t}w)\Big{|}^p\rmd w\right)^{\frac{1}{p}}\nonumber \\
&=\left(\int \rme^{-Nt}(1+|\rme^{-t}\chi|^m)^p\Big{|}f(\chi)\Big{|}^p\rmd w\right)^{\frac{1}{p}}\nonumber\\
&\leq\rme^{-\frac{Nt}{p}}\left(\int
(1+|\chi|^m)^p\Big{|}f(\chi)\Big{|}^p\rmd w\right)^{\frac{1}{p}} .
\label{eq:fetwweightedpnorm}
\end{align}

Putting \eqref{eq:phiweightlpnorm} together with
\eqref{eq:fetwweightedpnorm}, we can use Young's inequality in
\eqref{eq:dNucFPf} as before, since $1\leq q\leq p$ with $r$ given
by $\frac{1}{p}+1=\frac{1}{r}+\frac{1}{q}$ to get the desired
bound.
\endproof

\

\noindent {\bf Acknowledgements.-}
JAC and JR acknowledge partial support from DGI-MEC (Spain) project
MTM2005-08024 and 2005SGR00611 from Generalitat de Catalunya. We thank the Centre de Recerca Matem\`atica
(Barcelona) for partial funding and for providing an excellent
atmosphere for research.

\end{document}